# FINITE GROUPS WITH SMALL AUTOMIZERS FOR EVERY ABELIAN SUBGROUP OF NON PRIME POWER ORDER

RITESH DWIVEDI

Abstract. In this paper, we discuss about finite groups in which, $C_G H = N_G H$, for every abelian subgroup $H$ of non prime power order. Also, we classify all such nilpotent and minimal non nilpotent groups

## 1. Introduction

Throughout this paper, $G$ denotes a finite group and $p$ and $q$ denote a prime. By a $p$-group $G$, we mean $G$ is a group of prime power order, for some prime $p$. By a non $p$-group $G$, we mean $G$ is a group of non prime power order. Also, $C_G H$ and $N_G H$ denote the centralizer and the normalizer of $H$ in $G$, respectively. Also, we use $F(G)$ and $Z(G)$ to denote the Fitting subgroup and the center of $G$, respectively. Further, $[H]K$ denotes a split extension of $K$ by a normal subgroup $H$.

Recall that, the automizer of a subgroup $H$ in $G$ is defined as $Aut_G H = N_G H / C_G H$. Therefore $Aut_G H$ can be considered as a subgroup of $Aut H$ and since $Aut_G H$ contains an isomorphic copy of $Inn H$, we have $Inn H \leq Aut_G H \leq Aut H$. Now $Aut_G H$ is said to be small if $Aut_G H \cong Inn H$ and $Aut_G H$ is said to be large if $Aut_G H \cong Aut H$.

All finite groups with small automizers for their nonabelian subgroups are classified in [3]. Zassenhaus proved in [10] that $N_G H = C_G H$ for every abelian subgroup iff $G$ is abelian (also see [6]). Therefore $Aut_G H$ is small for all abelian subgroups of $G$ iff $G$ is abelian. Further, Li [8] classified all finite groups in which any nonmaximal abelian subgroups do not have any nontrivial inner automorphisms. Such groups are called $NC$-groups. All

---

finite groups, in which every non normal abelian subgroup has trivial automizer, are classified in [1]. Such groups are called quasi-$NC$ group. Finite groups, in which for any non normal abelian subgroup $A$, either $Aut_G A = 1$ or $C_G A = A$, are called $NNC$-groups and are discussed in [2]. More generally, all finite groups in which $Aut_G H$ is either small or large for every abelian subgroup H, are discussed in [9].

The following result is taken from [8, Lemma 2.4.]:

**Lemma 1.1.** *Let G be a group. Then the following are equivalent:*
*(1) G is an $NC$-group;*
*(2) $C_G A = A$ or $C_G A = N_G A$ holds for all abelian subgroups A of G of prime power order.*

The above lemma asserts also that if every abelian subgroup of $G$ of prime power order has small automizer, then infact every abelian subgroup of $G$ of non prime power order has also small automizer and $G$ become abelian. However, if every abelian subgroup of $G$ of non prime power order has small automizer, an abelian $p$-subgroup of $G$ need not have small automizer. For example, in the group $S_3 \times Z_3$, the only abelian subgroups of non prime power order are $\{I, \sigma\} \times Z_3$, where $\sigma$ is a transposition in $S_3$. Clearly all such subgroups have small automizer in $S_3 \times Z_3$. However, $A_3 \times Z_3$ is an abelian subgroup of prime power order but it's automizer is not small. Therefore it is of natural interest to classify all finite groups in which, for every abelian subgroup $H$, either $H$ is a $p$-subgroup ($p$ depends on $H$) or $C_G H = N_G H$. We shall do this partially in the present paper. We call such groups as $PNC$-group. Clearly every finite abelian group is $PNC$. So we shall focus on finite non abelian groups, which are $PNC$. The aim of this paper is to study all finite solvable $PNC$-groups. Also, we give a characterization of nilpotent $PNC$-groups and minimal non nilpotent $PNC$-groups.

## 2. Nilpotent　　$PNC$-gr oups

**Lemma 2.1.** *The class of $PNC$-groups is subgroup closed.*

*Proof.* It is obvious.　　　　　　　　　　　　　　　　　　　　　?

**Prop osition 2.2.** *Let G be a decomposable group with $G \cong G_1 \times G_2$, for some finite non trivial groups $G_1$ and $G_2$. If, either $(|G_1|, |G_2|) = 1$ or both*

$G_1$ and $G_2$ are non p-groups, then the following are equivalent:
(i) $G$ is abelian;
(ii) $G$ is $PNC$.

*Proof.* Clearly, if $G$ is abelian, then it is $PNC$. Conversely, assume that $G$ is $PNC$. First we observe that, under the given conditions, $N_{G_i}H_i = C_{G_i}H_i$, $\forall\ i = 1, 2$ and for every abelian subgroup $H_i$ of $G_i$. Now the result follows from [10, Theorem 7].                                     ?

**Corollary 2.3.** *Let $G$ be a non abelian group. Then $G$ is a nilpotent $PNC$-group iff $G$ is a p-group.*

*Proof.* If possible, suppose that $G$ is a non p-group. Then we have $G \cong G_1 \times G_2$, where $G_1$ is a Sylow-p subgroup of $G$ and $G_2$ is a Hall-$p^0$ subgroup of $G$. Also, we have $(|G_1|, |G_2|) = 1$. Now by above Proposition, $G$ is an abelian group, a contradiction. Hence $G$ is a p-group.
Converse is obvious.                                     ?

## 3. Sol vable  non  nilpotent    $PNC$-gr oups

Recall that the Fitting subgroup of $G$, denoted as $F(G)$, is defined as the largest normal nilpotent subgroup of $G$. Also, we denote the commutator subgroup of $G$ as $G^0$.

**Lemma 3.1.** *Let $G$ be a solvable non abelian group. If $G$ is a $PNC$-group, then $F(G)$ is a p-group.*

*Proof.* Note that $F(G)$ is a nilpotent $PNC$-group. Now if $F(G)$ is non abelian, then by Corollary 2.3, it is a p-group. Also, if $F(G)$ is an abelian but non p-group, then $C_G F(G) = N_G F(G)$. Since $G$ is solvable, we have $F(G) = C_G F(G) = N_G F(G)$. This gives that $G$ is abelian, a contradiction. This completes the proof.                                     ?

**Lemma 3.2.** *Let $G$ be a solvable non abelian group. If $G$ is a $PNC$-group, then $Z(G)$ is a p-group.*

*Proof.* First we note that $Z(G) \subseteq F(G)$. Now the result follows from Lemma 3.1.                                     ?

**Lemma 3.3.** *Let $G$ be a non abelian group with $G^0$ nilpotent. If $G$ is a $PNC$-group, then $G^0$ is a p-group.*

*Proof.* First we note that if $G^0$ is nilpotent then $G^0 \subseteq F(G)$. Now the result follows from Lemma 3.1. □

Recall that a group is said to be supersolvable, if it has a normal series with cyclic factors.

If $G$ is $PNC$, $F(G)$ need not be a Sylow subgroup. For example, $S_4$ is $PNC$ but it's Fitting subgroup is not a Sylow subgroup. However we have the following result.

**Lemma 3.4.** *Let $G$ be a non abelian supersolvable group. If $G$ is $PNC$, then $F(G)$ is a Sylow-p subgroup.*

*Proof.* First we note that the commutator subgroup of a supersolvable group is nilpotent. Now by above Lemma, $G^0$ is a $p$-group and hence $F(G)$ is a Sylow-$p$ subgroup. □

**Lemma 3.5.** *Let $G$ be a $PNC$-group and $H \leq Z(G)$. If $H \cap G^0 = 1$, then $G/H$ is also $PNC$.*

*Proof.* Let $K/H$ be an abelian non $p$-subgroup of $G/H$. Then $K$ is also a non $p$-subgroup of $G$. Also, $K^0 \subseteq H \cap G^0 = 1$. This gives that $K$ is abelian. Therefore $N_G K = C_G K$. Now, $N_{G/H}(K/H) = N_G K/H = C_G K/H \leq C_{G/H}(K/H)$. This implies that $N_{G/H}(K/H) = C_{G/H}(K/H)$. □

Recall the following results from [1]:

**Theorem 3.6.** *(see[1], Theorem 4.1.) Let $G$ be a non-nilpotent quasi-$NC$ group. Then all Sylow subgroups of $G$ are abelian.*

**Theorem 3.7.** *(see[1], Theorem 4.4.) Suppose that $G$ is non nilpotent quasi-$NC$ group. Then $G/Z(G)$ is also a quasi-$NC$ group.*

Note that if $G$ is a non nilpotent $PNC$-group, then Sylow subgroups of $G$ need not be abelian. For example, $S_4$ is a non nilpotent $PNC$-group, but it's Sylow-2 subgroups are non abelian. However the following result is an analogous of Theorem 3.7. Recall that, an $A$-group is a finite group with the property that all of its Sylow subgroups are abelian.

**Lemma 3.8.** *Let $G$ be a non nilpotent $A$-group. If $G$ is $PNC$, then $G/Z(G)$ is also $PNC$.*

*Proof.* From [1, Lemma 2.7], $G^0 \cap Z(G) = 1$. Now the result follows from Lemma 3.5. □

Recall that a finite group $G$ is said to be $SBP$, if every proper subgroup $H$ of $G$ has prime power order (the prime depending on $H$). The following classification of $SBP$ groups is given in [4; chapter 3, page 74/75].

**Theorem 3.9.** *Let $G \neq 1$ be an $SBP$ group. Then precisely one of the following holds:*
*(I) $G$ is a p-group, for some prime p.*
*(II) $|G| = pq$, where p and q are distinct primes.*
*(III) $|G| = p^a q$, p and q are distinct primes, $a = exp(p,q)$ and $a \geq 2$. More ever, $G \cong Z_p^a \times_\theta Z_q$ for some monomorphism $\theta : Z_q \to Aut Z_p^a$.*

Recall that a group $H$ is said to a central extension of $G$, if $H/A \cong G$, where $A \subseteq Z(H)$. For our convenience, we call a central extension $H$ of $G$ as a $p$-central extension, if $Z(H)$ is a $p$-group. Clearly, an $SBP$-group is a solvable $PNC$-group and if it is non abelian, then it's Fitting subgroup is a $p$-group. Now we prove the following result.

**Lemma 3.10.** *Let $G$ be an $SBP$-group with $F(G)$ a p-group. Then a p-central extension of $G$ is $PNC$.*

*Proof.* If $G$ is an abelian group or a $p$-group, the result holds trivially. So let $G$ be a non abelian group of order $p^n q$, with $F(G)$ a $p$-group. Also let $H$ be a $p$-central extension of $G$ with $H/A \cong G$, where $A \subseteq Z(H)$. Since $G$ is $SBP$, $A = Z(H)$. Also, since we have $F(H/Z(H)) = F(H)/Z(H)$, $F(H)$ is a $p$-group. Now we claim that every proper non $p$-subgroup of $H$ is of the form $TQ$, where $T \subseteq Z(H)$ and $Q$ is a Sylow-$q$ subgroup of $H$. Let $K$ be a non $p$-subgroup of $H$. Then $KA/A$ is a subgroup of $H/A$ and since $H/A$ is $SBP$, $|K/K \cap A| = |KA/A| = q$. This gives that $K = TQ$, where $T = K \cap A$ and $Q$ is a Sylow-$q$ subgroup of $K$. Therefore, for given any abelian non $p$-subgroup $TQ$ of $H$, $N_H(TQ) = C_H(TQ) = AQ$. This completes the proof. ?

Recall that a group is called a $CP$-group if every element of the group has prime power order (see [5]). Finite $CP$-groups were first studied by Higman [7] in 1957. The finite soluble $CP$-groups were classified by him as follows (see [7], Theorem 1):

**Theorem 3.11.** *Let $G$ be a soluble group all of whose elements have prime power order. Let p be the prime such that $G$ has a normal p-subgroup greater than 1, and let $P$ be the greatest normal p-subgroup of $G$. Then $G/P$ is either*

*(i) a cyclic group whose order is a power of a prime other than p; or (ii) a generalized quaternion group, p being odd; or (iii) a group of order $p^a q^b$ with cyclic Sylow subgroups, q being a prime of the form $kp^a + 1$. Thus G has order divisible by at most two primes, and G/P is metabelian.*

Now we prove some results about $CP$-groups.

**Lemma 3.12.** *Every finite CP-group is PNC.*

*Proof.* Note that a finite group is $CP$ iff every abelian subgroup of $G$ is of prime power order. Now the result follows from definition.                    ?

**Remark.** Unlike the case of $SBP$-groups, if $G$ is solvable $CP$-group with $F(G)$ a $p$-group, a $p$-central extension of $G$ need not be $PNC$. For example, $S_4$ is a solvable $CP$-group but $S_4 \times Z_2$ is not a $PNC$-group as $A_3 \times Z_2$ has not small automizer. However we have the following result.

**Theorem 3.13.** *Let G be a solvable CP-group of order $p^n q$ with $p > q$. Then a p-central extension of G is a PNC-group.*

*Proof.* Let $H$ be a $p$-central extension of $G$, with $H/A \cong G$, where $A \subseteq Z(H)$. Since $G$ is $CP$, $A$ can not be a proper subgroup of $Z(H)$ and so $A = Z(H)$. Let $K$ be an abelian non $p$-subgroup of $H$. Then $KA/A$ is an abelian subgroup of $H/A$ which is $CP$, being isomorphic to $G$. Therefore $|K/K \cap A| = |KA/A| = q$. This gives that $K = TQ$, where $T = K \cap A$ and $Q$ is a Sylow-$q$ subgroup of $K$. Clearly, $AQ \subseteq C_H(TQ)$. Now we shall show that $N_H(TQ) \subseteq AQ$. First note that $N_G S = S$, where $S$ is a Sylow-$q$ subgroup of $G$. For if $|N_G S| = p^r q$, for some positive integer $r \geq 1$. Then, since $p > q$, $G$ will contain an element of order $pq$, which is not possible. Hence $AQ/A = N_{H/A}(AQ/A) = N_H(AQ)/A$. This implies that $N_H(AQ) = AQ$. Now we have $N_H(TQ) \subseteq N_H(AQ) = AQ$. Hence $N_H(TQ) \subseteq AQ \subseteq C_H(TQ)$. This completes the proof.                    ?

**Prop osition 3.14.** *Let G be a non abelian solvable group with an abelian maximal normal subgroup. If G is a PNC-group, then $|G| = p^n q$, where p, q are primes (need not be distinct) and n is some positive integer.*

*Proof.* First note that a maximal normal subgroup of a solvable group is of prime index. Now if $H$ is an abelian maximal normal subgroup but not a $p$-group, then $N_G H = C_G H$. It gives that $G$ is abelian, a contradiction. Hence $H$ is a $p$-group. Now the result follows.                    ?

**Corollary 3.15.** $D_{2n}$ *is a PNC-group iff* $n = p^k$ *, for some prime p and positive integer k.*

*Proof.* Since $<r>$ (where $r \in D_{2n}$ is an element of order $n$) is an abelian maximal normal subgroup of $D_{2n}$, from the above proposition $<r>$ is a $p$-group. Now the result follows. Conversely, it is easy to observe that, if $n = p^r$, then $D_{2n}$ is $PNC$. In fact in this case $D_{2n}$ is $CP$.                        ?

We have seen in the section 2 that nilpotent $PNC$ groups are $p$-group. However if we impose a weaker condition on $G$, namely the condition of minimal non nilpotent, then we have the following result.

**Theorem 3.16.** *G is a minimal non nilpotent PNC-group iff* $G \cong [F(G)]Q$, *where $F(G)$ is a Sylow-p subgroup of G and Q is a cyclic subgroup of order q acting non trivially on $[F(G)]$ and every proper non p-subgroup of G is abelian.*

*Proof.* Let $G$ be a minimal non nilpotent $PNC$ group. Then $G$ is solvable and so $F(G)$ is a $p$-group. Since every proper subgroup is nilpotent, $F(G)$ is a maximal normal subgroup of $G$ and so $[G : F(G)] = q$, for some prime $q \ne p$. Also from Corollary 2.3, it follows that every proper non $p$-subgroup of $G$ is abelian. Since $G$ is non nilpotent, Sylow-$q$ subgroup is non normal. So let $Q$ be a Sylow-$q$ subgroup of $G$. Then clearly, $G \cong [F(G)]Q$ and $Q$, being non normal, acts non trivially on F(G).
Conversely, assume that $G \cong [F(G)]Q$, where $F(G)$ is a Sylow-$p$ subgroup of $G$ and $Q$ is a cyclic subgroup of order $q$ acting non trivially on $[F(G)]$ and every proper non $p$-subgroup of $G$ is abelian. Since the action of $Q$ on $F(G)$ is non trivial, $G$ is non nilpotent. Also every non $p$-subgroup of $G$ is abelian, so $G$ is minimal non nilpotent. Now since $F(G)$ is a $p$-group, no proper non $p$-subgroup of $G$ is normal. Therefore $N_G(H) = C_G(H)$, for every abelian non $p$-subgroup of $G$. This completes the proof.            ?

**Acknowledgemen t:** I am grateful to the referee for the comments which improved the quality of the paper.

Ritesh Dwivedi

Prayagraj, India.

E-mail: riteshgrouptheory@gmail.com